\documentclass[11pt]{article}
\usepackage{amssymb}
\usepackage[dvips]{graphics}
\usepackage{theorem}

\newtheorem{theorem}{{\sc Theorem}}[section]

\newtheorem{proposition}[theorem]{{\sc Proposition}}
\newtheorem{conjecture}{{\sc Conjecture}}
\newtheorem{definition}{{\sc Definition}}[section]
\newtheorem{lemma}{{\sc Lemma}}[section]
\newtheorem{corollary}[theorem]{{\sc Corollary}}
\newtheorem{remark}[theorem]{{\sc Remark}}
\newtheorem{qquestion}[theorem]{{\sc Question}}
\newtheorem{pprobleme}[theorem]{{\sc Problem}}
\newtheorem{notation}{{\sc Notation}}

\newenvironment{question}{\noindent {\sc Question} \begin{em}}
{\end{em}\nolinebreak }

\newenvironment{proof}{{\sc Proof}}{\nolinebreak $\Box $}

\newenvironment{demonstration}{{\sc Demonstration}} {{\sc QED}}

\newenvironment{resume}{\small \begin{center} {\bf R\'esum\'e}
\end{center} \hspace{1cm} \begin{minipage}[t]{13cm} \hspace{.5cm} } 
{\end{minipage} \normalsize }

\begin{document}

\title{Two transformations of simple polytopes preserving moment-angle 
manifolds}
\author{Fr\'ed\'eric Bosio}
\date{\today}
\maketitle

\begin{resume}
  Nous introduisons ici deux constructions assez g\'en\'erales sur les 
polytopes simples, conduisant \`a des polytopes ayant m\^emes vari\'et\'es 
moment-angle. Comme application, nous donnons des examples de polytopes (dual) 
flag qui ne sont pas rigides.
\end{resume}

\begin{abstract}
  In this paper, we introduce two new, quite frequently occuring, 
constructions, yielding polytopes with diffeomorphic moment-angle manifolds. As 
an application, we give examples of (dual) flag polytopes that are not rigid.
\end{abstract}

\section*{Introduction}

  Moment-angle manifolds form a central object in toric geometry. They have 
been introduced by Buchstaber and Panov \cite{B-P} to oversee the so-called 
quasitoric manifolds. To any simple (combinatorial) polytope is associated such 
a manifold that supports a toric action whose orbit space is the given polytope.

  The geometry of the moment-angle manifold associated to a polytope is then 
completely determined by the combinatorics of the polytope. For example, we can 
describe the homology of a moment-angle manifold in terms of the polytope 
(see \cite{Ba}) and we have more precise results in some cases (\cite{L-V},
\cite{B-P}, \cite{B-M}, \cite{L-G}, \cite{Ir}...)

  A general and very intriguing problem can (roughly) be stated as follows:

\begin{qquestion}
When do two polytopes produce "the same" moment-angle manifolds?
\end{qquestion}

  Closely related to this question is the following one:

\begin{qquestion}
\label{nonrigid}
  Given a polytope, are there other polytopes giving the same moment-angle 
manifold?
\end{qquestion}

  A polytope can be thought of as {\em{rigid}} if the answer to 
question~\ref{nonrigid} is negative.

\vskip 3mm

  Different quite natural ways of being "the same", more or less strong, can be 
considered.

  Recall that any moment-angle manifold supports the natural action of a torus 
whose orbit space is the polytope in question.

  Given two simple polytopes $P$ and $Q$, we know that $Z_P $ and $Z_Q $ are 
equivariantly diffeomorphic (with respect to the forementioned actions) if and 
only if $P$ and $Q$ are combinatorially equivalent \cite{B-M}. Hence, this 
request is too strong to involve an interesting problem.

  We can weaken the request by only demanding that the moment-angle manifolds 
are diffeomorphic or even that their cohomology rings are isomorphic (with the 
natural graduation of the cohomology of a manifold). The first request seems 
much stronger than the second, but the equivalence between them (and so between 
other notions of intermadiate strength) is an open question.

  We use here the notion of graded diffeomorphy. This notion is stronger than 
diffeomorphy, in which some important discrete invariants may vary (see 
\cite{Bo1}). Nevertheless, the most widespread constructions preserving the 
differential structure of moment-angle manifolds also preserve the bigraduation.

  In this paper, we introduce two new, quite frequently occuring, 
constructions, yielding graded-equivalent polytopes. The first one relies on 
the notion of flip of polytopes (which can be thought as "the smallest change 
that a polytope can undergo"). It is called a biflip as it is a succesions of 
two flips (not any). The second one is called "puzzle-move" as it relies on 
an assembly of pieces. As an application, we give the first examples of 
nonrigid (dual) flag polytopes.

\section{Preliminaries}

\paragraph{Polytope}

  We only consider here combinatorial polytopes. We consider a simple 
$d$-dimensional polytope $P$ as the set of all subsets of facets of a geometric 
simple $d$-dimensional polytope $P_{geom}$ whose intersection are nonempty.

  In this sense:

- The elements of the elements of $P$ are identified with the facets of 
$P_{geom}$ also called the facets of $P$. They are called the facets of $P$ and 
we identify a facet $F$ of $P$ with the element $\{ F \} $ of $P$. We usually 
note ${\cal F }$ the set of facets of $P$.

- The maximal elements of $P$ have $d$ elements and correspond to the subsets 
of facets of $P_{geom}$ containing a given vertex. They are called the vertices 
of $P$.

- An element of $P$ is called a face of $P$ (and corresponds to the subsets of 
facets of $P_{geom}$ containing a given face).

  Recall just that a combinatorial polytope is completely determined by its 
vertices.

\paragraph{Moment-angle manifold}

\begin{definition}
  Let ${\cal E }$ a finite set. We note $T^{\cal E }$ the real torus 
$S^{\cal E }$ where $S$ denotes the unit circle in $\mathbb C $, with its usual 
structure of topological group.
\end{definition}

Let's begin by recall a construction of a moment-angle manifold. We start 
from a simple polytope $P$, and denote ${\cal F }$ the set of its facets. We 
can consider the torus $T^{\cal F }$ and the product $P \times T^{\cal F }$, on 
which $T^{\cal F }$ acts by translation on itself.

  On $P \times T^{\cal F }$, we identify two elements $(p,t)$ and $(p',t')$ if 
$p=p'$ and $p$ lies on every facet not sent on $1$ by $t^{-1}t'$.

  The space $Z_P $ obtained by this identification is the moment-angle 
manifold over $P$. The action of $T^{\cal E }$ descends on $Z_P $. (We can keep 
in mind that each facet correspond to a rotation of $Z_P $ fixing the 
submanifold over this facet, and that all these rotations commute).

\paragraph{Homology}

  We recall (see \cite{Ba}) that the cohomology of $Z_P $ is endowed with a 
natural bigraduation, more precisely, for any subset ${\cal X }$ of facets of 
$P$, the homology of the space $\bigcup_{F \in {\cal X }} F$ corresponds to a 
part (a submodule) of the homology of $Z_P $, and we recall the following 
formula~:
$$H_k (Z_P , {\mathbb Z }) \simeq 
\bigoplus_{{\cal X } \subset {\cal F}} 
\tilde{H}_{k - |{\cal X }| - 1} (P_{\cal X } , {\mathbb Z })$$

  We note $H^{p,q} (Z_P , \mathbb Z )$ the subspaces of the homology of $Z_P $ 
on whose the bigraduation has a given value.

\begin{definition}
  A graded diffeomorphism $\phi $ between $Z_P $ and $Z_{P'}$ is a 
diffeomorphism between them such that, for any $p,q$, we have 
$\phi ^* (H^{p,q} (Z_{P'}, \mathbb Z )) =  H^{p,q} (Z_P , \mathbb Z )$.

  Two polytopes $P$ and $Q$ are said $Gr$-equivalent if there is a graded 
diffeomorphism between their moment-angle manifolds.
\end{definition}

\paragraph{Flip}

\begin{definition}
  Consider a simple polytope $P$ of which $A$ is a simplicial face.

  A facet $F$ of $P$ is called a {\em{bounding facet}} of $A$ if it meets the 
boundary of $A$ but not its (relative) interior.
\end{definition}

  Here, we only consider flips of polytopes at combinatorial level. Consider 
a polytope $P$ of dimension $d \geq 2$. Let $1 \leq p \leq d$ an integer, and
$A$ a simplicial $p-1$ face of $P$. Now consider the set of facets 
${\cal A } = \{ A_1 ,..., A_{d-p+1} \} $ so that 
$A = \bigcap_{F \in {\cal A }} F$ and the set 
${\cal B } = \{ B_1 ,..., B_p \} $ of bounding facets of $A$. Call now $V$ the 
set of vertices of the simplex on ${\cal A } \cup {\cal B }$ (i.e. the set of 
complements of singletons in ${\cal A } \cup {\cal B }$).

  Now, consider the symmetric difference between the set of vertices of $P$ and 
$V$. If it is the list of vertices of a polytope $Q$, we say that $Q$ is 
obtained from $P$ by a flip, and setting $q = d-p+1$, we say this flip is a 
$(p,q)$-flip.

  If $Q$ is obtained from $P$ by a $(p,q)$-flip, then $P$ is obtained from $Q$ 
by a $(q,p)$-flip where the roles of ${\cal A }$ and ${\cal B }$ are inverted.

\includegraphics{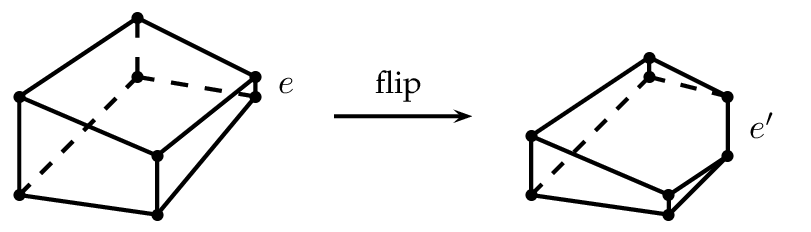}

\paragraph{Wedge}

  We give here a purely combinatorial presentation of the wedging operation. It 
can be also presented more geometrically, for istance in~\cite{K-W}.

  Let $P$ a simple $d$-polytope and $F$ a facet of $P$. We note $W_F P$, called 
the {\em{wedge over $P$ on $F$ }}, the simple $d+1$-polytope obtained by 
replacing $F$ by two elements $F_1 $ and $F_2 $, and whose vertices are:

- the sets $\{ F_1 ; G_1 ;...; G_d \} $ or 
$\{ F_2 ; G_1 ;...; G_d  \} $ where $\{ G_1 ;...; G_d \} $ is a vertex of $P$

- the sets $\{ F_1 ; F_2 ; G_1 ;...; G_{d-1} \} $ where 
$\{ F ; G_1 ;...; G_{d-1} \} $ is a vertex of $P$.

  The facets $F_1 $ and $F_2 $ are called the {\em{main facets}} of the wedge.

  This operation can be repeated as many times as we want over any facet, and 
to a map $\alpha : {\cal F } \to {\mathbb N }$, we note $W_{\alpha } P$ the 
polytope called multiwedge over $P$ the polytope obtained by taking 
$\alpha (F)$ times a wedge on $F$ for each $F$, which is well defined, i.e. 
does not depend on the order in which wedges are performed.

\section{Biflips}

  Consider a simple polytope $P$. We know that two polytopes obtained by 
different flips of equal index from $P$ may not be $Gr$-equivalent. For 
instance, if we consider the pentagonal book, whose associated moment angle 
manifold is diffeomorphic to 
${\# \atop 3} S^3 \times S^6 {\# \atop 2} S^4 \times S^5 $, an edge flip may 
produce either another pentagonal book or a cube, whose associated moment angle 
manifold is $S^3 \times S^3 \times S^3 $.

  There is nevertheless a simple case in which it is possible to guarantee that 
two flips will produce $Gr$-equivalent polytopes:

\begin{theorem}
  Let $P$ a simple polytope. Consider two simplicial faces $F_1 $ and $F_2 $ 
and assume that their bounding facets are the same.

  Then if both faces are flippable, the two polytopes obtained after the flips 
are $Gr$-equivalent.
\end{theorem}

  This theorem leads to the notion of biflip:

\begin{definition}
  Let $Q$ a polytope. A sequence of two flips 
$Q \stackrel{f_1 }{\rightarrow} P \stackrel{f_2}{\rightarrow} R$ is called a 
biflip if the bounding facets of the simplicial face of $P$ appearing by $f_1 $ 
are the same as the bounding facets of the one disappearing by $f_2 $.
\end{definition}

\rotatebox{90}{\scalebox{.45}{\includegraphics{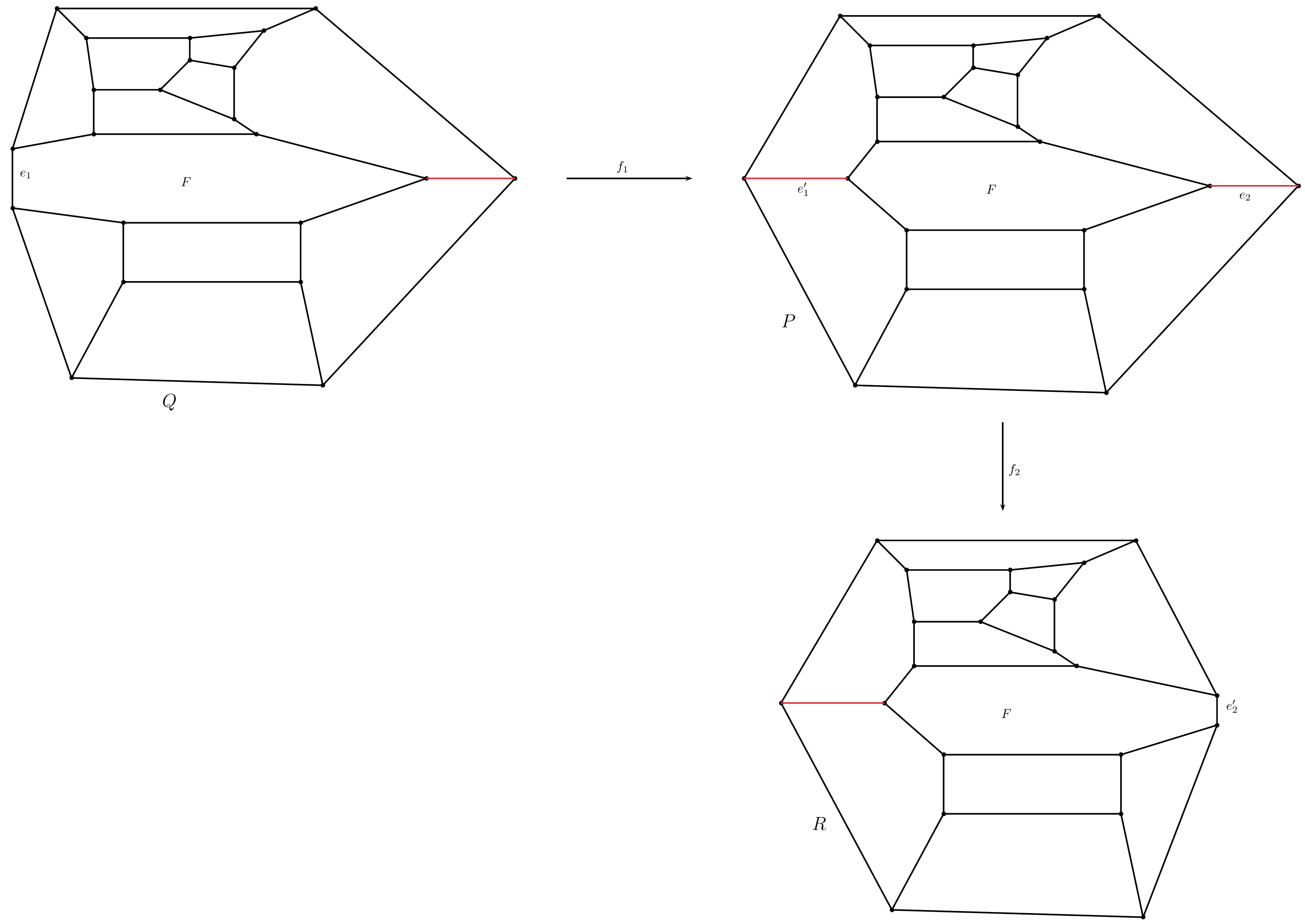}}}

  We immediately deduce from the theorem:

\begin{corollary}
  Two polytopes joined by a (sequence of) biflip(s) are $Gr$-equivalent.
\end{corollary}

  Let's prove the theorem.

\begin{demonstration}

  As explained in \cite{B-M}, it is possible to describe the change on the 
geometry of a moment angle manifold arising after a flip. This is given by a 
equivariant surgery. In the forecited article, it is shown that any 
$(1,n)$-flip on a polytope induces the same (non fully equivariant) surgery, 
hence the same differential structure. The proof is based on the fact that some 
submanifolds are isotopic inside the moment angle manifold.

  In the case of more general flips of the same index, the topology of the sets 
that are removed and glued, as well as the glueing operation, only depend on 
the index. In this sense, there only remains to prove that the removed pieces 
are isotopic inside the moment-angle manifold $Z_P $.

  If $F_1 $ and $F_2 $ are $(k-1)$-simplices, these removed pieces have the 
form $S^{2k-1} \times T^{n-d+1} \times D^{2(d-k)}$, indeed are a product 
neighbourhood of the submanifold $S^{2k-1} \times T^{n-d+1}$ over $F_1 $ or 
$F_2 $.

  First, remark that a facet that does not bound a given simplicial face, 
either contains it or is disjoint from it.

  Consider then the bounding facets ${\cal A }$ of $F_1 $ and $F_2 $. Call also 
${\cal B }_0 $ the facets containing $F_1 $ and $F_2 $, ${\cal B }_1 $ those 
containing $F_1 $, not $F_2 $, ${\cal B }_2 $ those containing $F_2 $, not 
$F_1 $, and ${\cal C }$ those disjoint from both $F_1 $ and $F_2 $.

  Call $k = | {\cal A } |$, so $F_1 $ and $F_2 $ have dimension $k-1$. Call 
$r = | {\cal B } |$, so $| {\cal B }_1 | = | {\cal B }_2 | = d-k-r+1$
and $| {\cal C } | = n-2d+k+r-2$.

  For a bounding facet $A$ in ${\cal A }$, note $v_{A,1}$ (resp. $v_{A,2}$) the 
vertex of $F_1 $ (resp. of $F_2 $) not lying on $A$.

  Consider the simplex $\Delta ^{k-1}$ on $k$ vertices, each vertex $e_A $ 
being associated to a facet $A$ in ${\cal A }$.

  We construct now an explicit isotopy between the forementioned loci. We 
first consider the map $g$ from $\Delta ^{k-1} \times [0;1]$ to $P$ that is 
separately affine and so that for all $A$, $g(e_A ,0) = v_A $ and 
$g(e_A ,1) = v'_A $. Notice in particular that $\Delta ^{k-1} \times \{ 0 \} $ 
is sent on $F_1 $ and $\Delta ^{k-1} \times \{ 1 \} $ on $F_2 $.

  We fix a bijection $c$ between ${\cal C }$ and $[1 ,..., n-2d+k+r-2]$ as
well as two bijections $b_i ,i=1,2$ between respectively ${\cal B }_i $ and 
$[1 ,..., d-k-r+1]$, then consider the map $\phi $ from the torus $T^{n-d-1}$ 
to $T^{\cal F }$ given by 
$\phi (t_1 ,..., t_{d+1-k}) = (x_F )_{F \in {\cal F }}$ where $x_F = 1$ if $F$ 
is in ${\cal A }$, $x_F = t_{b(F)}$ if $F$ is in ${\cal B }$, 
$x_F = t_{r + b_i (F)}$ if $F$ is in ${\cal B }_i , i=1,2$ and 
$x_F = t_{(d-k-r+1) + c(F)}$ if $F$ is in ${\cal C }$.

  Notice $\phi $ is a group morphism and an embedding. Now, we define a map $h$ 
from $\Delta ^{k-1} \times T^{n-d+k-1} \times [0;1]$ to $Z_P $ by considering 
$T^{n-d+k-1}$ as $T^{\cal A } \times T^{n-d+1}$, and posing 
$h(x , \gamma _A , \gamma , t) = 
\left( \phi (\gamma )  \: \gamma _A \right) \cdot g(x,t)$.

  This action descends to $S^{2k-1} \times T^{n-d+1} \times [0;1]$. For $t=0$, 
this amounts to the action of $T^{{\cal B }_2 } \times T^{\cal C }$ on 
$S^{2k-1}$ giving the submanifold of $Z_P $ over $F_1 $ and for $t=1$, we get 
the action of $T^{{\cal B }_1 } \times T^{\cal C }$ on $S^{2k-1}$ giving the
submanifold of $Z_P $ over $F_2 $. Thus these two submanifolds are isotopic 
inside $Z_P $.

  Moreover, we immediately extend this isotopy to arbitrary small Collar 
neighbourhoods of these two submanifolds and easily entend this to a global 
isotopy of diffeomorphisms of $Z_P $, all being identity outside arbitrary 
small neighbourhoods of the image of $h$. This proves that $Z_{P_1 }$ and 
$Z_{P_2 }$ are diffeomorphic.

  We now verify that our diffeomorphism preserves the bigraduation.

  The homology of $Z_{P_1 }$ is generated by four kinds of classes~:

\begin{enumerate}

\item Images of classes that are represented by a cycle which is disjoint from 
$F'_1 $.

\item The top-class.

\item Images of classes for whose any representant is homologous to the 
boundary of $F'_1 $.

\item Images of (pure) classes whose intersection with some class of the 
precedent type equal $1$.

\end{enumerate}

  This first case corresponds to the case the representing cycle is inside the 
locus of $Z_{P_1 }$ where the diffeomorphism is identity. So the same cycle 
represents also the image class.

  The second case is immediate, as both polytopes have the same dimension.

  In the third case, The set of facets on which the representant is drawn must 
contain all facets bounding $F'_1 $, and no bounding facet of $F_1 $ otherwise 
it would be a boundary. Turn now in ${\cal X }$ these facets into the bounding 
facets of $F'_2 $ and the boundary of $F'_1 $ into the one of $F'_2 $. We 
preserve the degree and these classes correspond by the diffeomorphism.

  In the last case, a representant of such a class must be linked with the 
boundary of $F'_1 $ and a representant of its image must be linked with the 
boundary of $F'_2 $. So here too the degrees must fit.

  This completes the proof of the theorem.
\end{demonstration}

\begin{remark}
  Two polytopes obtained by cutting a vertex to the same polytope $P$ are 
joined by a biflip as cutting a vertex produces a $(1,d)$-flip. So we recover 
that two such polytopes are $Gr$-equivalent.
\end{remark}

  Let's give an application of this notion.

  Recall that the list of simplicial $4$-polytopes with $8$ vertices appears in 
\cite{G-S}, where they are numbered from $P_1 ^8 $ to $P_{37}^8 $. We know 
\cite{Ir} that some moment-angle manifold, viz. the one associated to the dual 
of the simplicial polytope $P_{28}^8 $, is a connected sum of sphere products, 
with a summand being the product of three spheres, whereas all (nontrivial) 
previously known examples of this kind had every summand being a product of two 
spheres.

\begin{proposition}
  The dual of polytope $P_{28}^8 $ is not rigid, more precisely it is 
$Gr$-equivalent to the duals of $P_{27}^8 $ and $P_{29}^8 $ (and to no other 
polytope).
\end{proposition}

\begin{proof}
  Indeed, we will show that the three polytopes in question are joined by 
biflips.
  
  We remark that these three polytopes are exactly the three simple polychora 
with four heptahedra and four hexahedra (whithout other facet).

  Consider a dual cyclic polychoron with seven facets. Number its facets from 
$1$ to $7$ in the natural cyclic order. Truncate now a triangular face, say 
$1 \cap 3$. We could see the dual of the polychoron we get is $P_{22}^8 $. The 
list of its vertices is the following:
$$(1245) (1248) (1256) (1267) (1278) (1457) (1478) (1567) (2345) (2348)$$
$$(2356) (2367) (2378) (3456) (3467) (3478) (4567)$$
We remark that three edges are bounded by facets $5$ and $8$, namely $(124)$ 
(i.e. $1 \cap 2 \cap 4$), $(147)$ and $(234)$.

  We easily see that the three polytopes obtained by flipping one of these 
edges have only $\{ 1;3 \} $ and $\{ 6;8 \} $ as pairs of disjoint facets, so 
these four facets are hexahedra, whereas the other four ones are heptahedra. 
Moreover no two of these polychora are isomorphic. Indeed, the last one has a 
cubic facet ($3$), opposed to both others, and the second one has to facets 
that are pentagonal prisms ($2$ and $7$), opposed to both others. Hence we get 
our three polytopes, that are, as announced, joined by biflips, hence 
$Gr$-equivalent.
\end{proof}

\section{Puzzle-equivalence}

  We introduce here a new construction leading to $Gr$-equivalent polytopes. 
This notion that we call puzzle-equivalence is based on a generalization of the 
notion of blending (connected sum) of polytopes that has been introduced by 
Holt \cite{Ho}.

\subsection{Wedge equivalence}
  
\begin{definition}
  Two facets $F$ and $G$ of a polytope $P$ are called wedge-equivalent if the 
transposition of these two facets is an automorphism of $P$.
\end{definition}

  Clearly, wedge-equivalence is an equivalence relation, whose classes are 
called wedge-classes.

\begin{proposition}
\label{equivwedgeq}
  Let $F$ ang $G$ two facets of a simple polytope $P$. The following assertions 
are equivalent:

\noindent
i) $F$ and $G$ are wedge-equivalent; \newline
ii) any vertex of $P$ lies either on $F$ or on $G$; \newline
iii) $F$ and $G$ belong to the same missing faces of $P$.
\end{proposition}

\begin{proof}
  Indeed, if the transposition swapping $F$ and $G$ is not an automorphism of 
$P$, then there are facets $H_i , 1 \leq i \leq d-1$ different from $G$ so that 
$F \cap H_1 \cap ... \cap H_{d-1}$ is a vertex and 
$G \cap H_1 \cap ... \cap H_{d-1 }$ not. A minimal subset $G$ of 
$\{ H_i \} _{1 \leq i \leq d-1}$ whose intersection with $G$ is empty yields a 
missing face containg $G$ and not $F$. So iii) $\Rightarrow $ i).

  Assume now $F$ and $G$ are wedge-equivalent. If we had a vertex neither on 
$F$ nor on $G$, we would have, by connexity of the graph of $P$, an edge 
joining such a vertex $v$ to a vertex $v'$ on $F$ or on $G$ (we can assume $F$, 
so $v'$ is not on $G$). The image $v''$ of $v'$ by the transposition 
$F \leftrightarrow G$ should also be adjacent to $v'$ and on $G$, not on $F$. 
Hence, $v$ and $v''$ would be two distinct vertices that are adjacent to $v'$ 
and not on $F$, which is impossible. So i) $\Rightarrow $ ii).

  Assume every vertex lies either on $F$ or on $G$, and let ${\cal E }$ a 
missing face containing $F$. Then there is an vertex on each facet of 
${\cal E } \setminus \{ F \} $. Such a vertex $v$ is on $G$ by hypothesis and 
let $v'$ adjacent to $v$ and not on $G$. Then $v'$ lies on $F$ and on every 
facet of ${\cal E } - \{ F,G \} $, so on every facet of ${\cal E } - \{ G \} $. 
As ${\cal E }$ is a missing face, ${\cal E } - \{ G \} $ is not ${\cal E }$, so 
$G$ belongs to ${\cal E }$. Hence ii) $\Rightarrow $ iii).

  This completes the proof of the proposition.
\end{proof}

\begin{remark}
  The two main facets of a wedge are wedge-equivalent. Also if $P$ is a 
polytope, the facets $P \times \{ 0 \} $ and $P \times \{ 1 \} $ of
$P \times [0;1]$ are wedge-equivalent. This in fact the only cases, and we 
could even consider the product with an interval as the wedge on a ghost facet.
\end{remark}

\begin{definition}
  Let $P$ be a polytope. We say that an automorphism $\phi $ of $P$ is 
{\em{harmless}} (or of wedge type) if it stabilizes wedge-classes, i.e. if, for 
any facet $F$ of $P$, $\phi (F)$ is wedge-equivalent to $F$.
\end{definition}

  In other words, a harmless automorphism is a composition of automorphisms 
that are transpositions.

  In fact, harmless automorphisms of polytopes are closely related to 
diffeomorphisms of moment-angle manifolds:

\begin{proposition}
\label{caract-harm}
  Let $\phi $ be an automorphism of a polytope $P$. The following assertions 
are equivalent:

\noindent
i) The automorphism $\phi $ is harmless; \newline
ii) the diffeomorphism $\phi $ induces on $Z_P $ is isotopic to identity; 
\newline
iii) the diffeomorphism $\phi $ induces identity on the homology of $Z_P $.
\end{proposition}

\begin{proof}
  Let's first prove i) $\Rightarrow $ ii). Consider a path of unitary 
automorphisms of ${\mathbb C }^2 $ from the identity to the automorphism 
permutating coordinates, for instance:
$$(z , w , \theta ) \to 
\frac{1}{2}((z+w) + e^{i \theta }(z-w) , (z+w) - e^{i \theta }(z-w)$$
  For $\theta = 0$, we get the identity. For $\theta = \pi $, coordinates are 
swapped.

  This path of isometries can be restricted to any sphere centered at the 
origin.

  Assume $\phi $ is an automorphism of a polytope $P$ which only transposes 
facets $F$ and $F'$. The path in $Diff (Z_P )$ preserving all coordinates but
those of $F$ and $F'$ and acting on this pair of coordinates like the path
thereabove clearly connects the identity to the automorphism of $Z_P $
permuting the desired coordinates. This proves that the automorphism of $Z_P $ 
induced by the permutation of $F$ and $F'$ is homotopic to identity.

Remark: We use here that $Z_P $ is filled by spheres on coordinates $F$ and 
$F'$, thanks to the model of intersection of quadrics with, in each quadric, 
the same coefficients for the two wedge-equivalent coodinates in question.

  The implication ii) $\Rightarrow $ iii) is a basic result of homotopy theory 
which is valid in a much broader context~\cite{Sp}.

  Let's prove iii) $\Rightarrow $ i). Assume we have an automorphism $\phi $
of a polytope $P$ whose induced diffeomorphism of $Z_P $ acts homologically 
trivially.

  Let ${\cal E }$ a missing face of $P$. This missing face induces a homology 
class of $Z_P $. This class must be fixed by $\phi _* $, so $\phi $ must 
preserve ${\cal E }$. Then $\phi $ is harmless by proposition~\ref{equivwedgeq}.

  This completes the proof of the proposition.
\end{proof}

\subsection{Puzzle-equivalence}

  We now introuce our manipulation on polytopes:

  We consider two polytopes $P_1 $ and $P_2 $ embedded in a $d$-dimensional 
euclidean space $E$. We consider a hyperplane $H$ of $E$, whose associated 
half-spaces will be noted $H_+ $ and $H_- $, intersecting $P_1 $ and $P_2 $ but 
not at any vertex. For $i=1,2$, the intersection $P_i \cap H$ is then a 
$d-1$-dimensional polytope $\Delta _i $, and $P_i $ is obtained by glueing the 
two parts $P_i \cap H_+ $ and $P_i \cap H_+ $ along their common facet 
$\Delta _i $.

  We assume there are isomorphisms $i_+ $ and $i_- $ between one one hand 
$P_1 \cap H_+ $ and $P_2 \cap H_+ $, on the other hand $P_1 \cap H_- $ and 
$P_2 \cap H_- $, both sending $\Delta _1 $ on $\Delta _2 $ (in particular 
$\Delta _1 $ and $\Delta _2 $ have to be isomorphic). We still note $i_+ $ and 
$i_- $ the restrictions of $i_+ $ and $i_- $ to $\Delta _1 $, and we can 
consider the automorphism $\phi = (i_- )^{-1} \circ i_+ $ of $\Delta _1 $.

  If $\phi $ is harmless, then we say that we pass from $P_1 $ to $P_2 $ by a 
puzzle-move. In this case, we pass from $P_2 $ to $P_1 $ by the inverse 
puzzle-move.

  In fact, the combinatorial data of $(P_1 , \Delta ^+ , \phi )$ completely 
determines the combinatorics of $P_2 $, where $\Delta ^+ $ is $\Delta _1 $ with 
the choice of the halfspace $H_+ $, and we sometimes define a puzzle-move with 
this data. In this case, noting $P$ for $P_1 $, $P_2 $ is noted 
$P_{\Delta ^+ , \phi }$.

  We also can notice that when $\phi $ is an involution, in particular a 
transposition, the choice of the halfspace has no importance and we can omit it.

\begin{definition}
  A puzzle-move will be called {\em{trivial}} or {\em{nontrivial}} according to 
whether $P_2 $ is isomorphic to $P_1 $ or not.

  The equivalence relation between polytopes generated by puzzle-moves is 
called {\em{puzzle-equivalence}}.
 \end{definition}

  Here is a picture of a puzzle-move:

\scalebox{.6}{\includegraphics{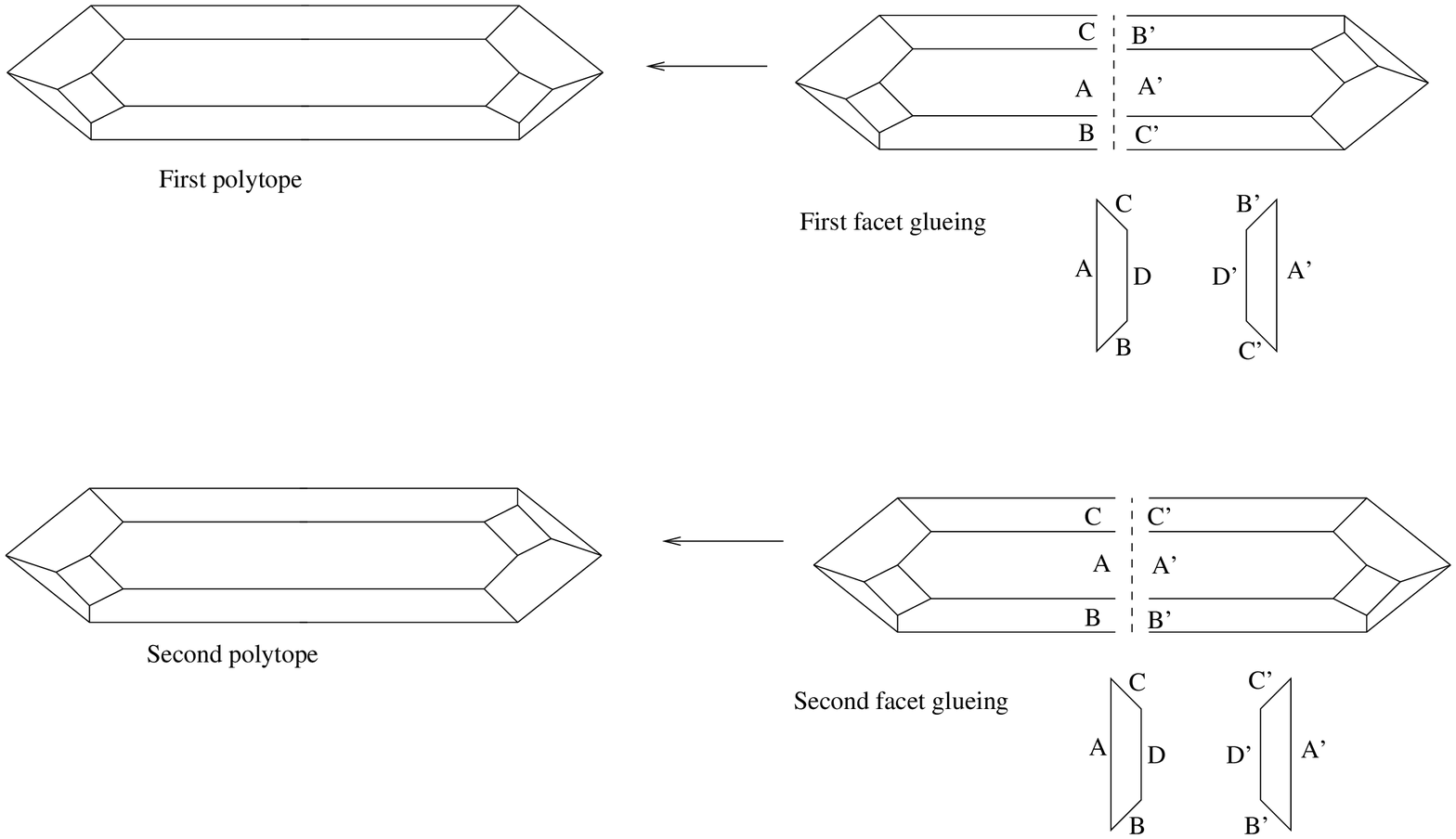}}

\begin{remark}
  We can notice that $P_{\Delta , \phi }$ actually depends on $\Delta $. We can 
construct examples of different puzzle-moves with the same $P$ and $\phi $.

\scalebox{.9}{\includegraphics{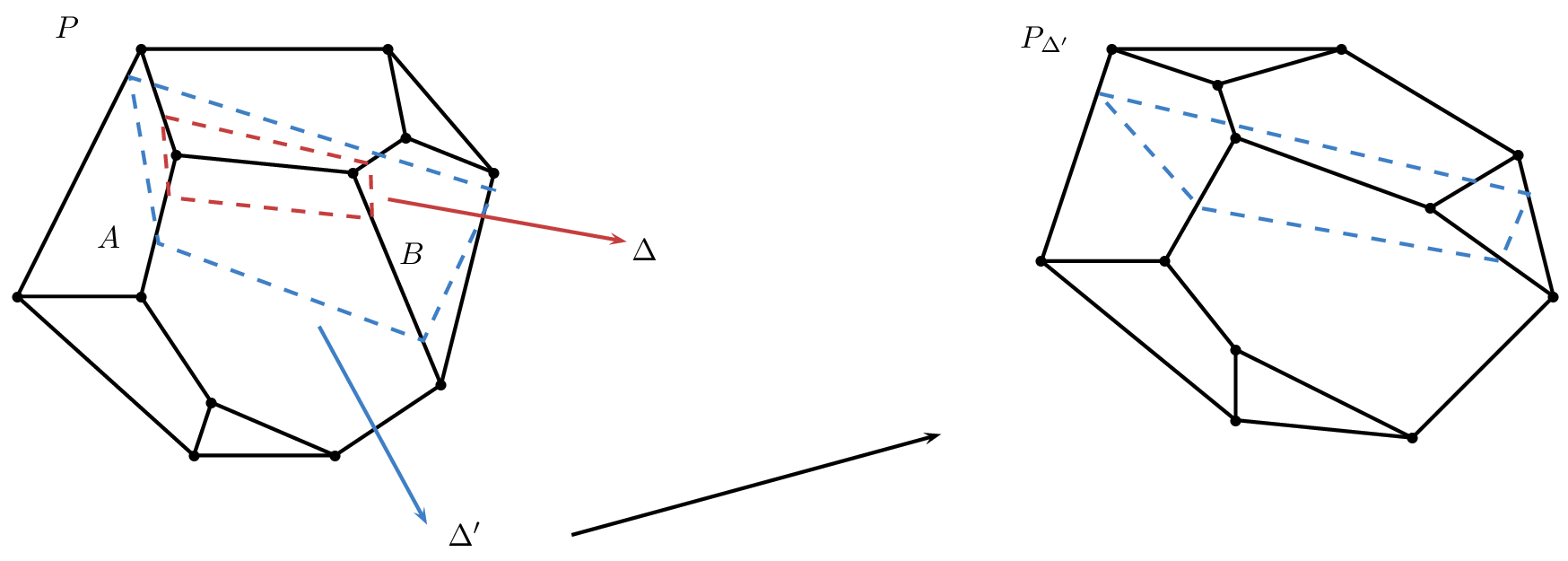}}

  The puzzle-move $(P , \Delta , (AB))$ is trivial whereas 
$(P , \Delta ' , (AB))$ is not.
\end{remark}

\begin{theorem}
  Two puzzle-equivalent polytopes ar $Gr$-equivalent.
\end{theorem}

\begin{remark}
  In particular, if $\Delta $ is a simplex, then any automorphism of $\Delta $ 
is harmless and any connected sum can be performed, so 
$(P , \Delta ^+ , \phi )$ is always a puzzle-move. Hence we recover that two 
polytopes obtained by connected sums of the same polytopes $P$ and $Q$ on the 
same vertices are $Gr$-equivalent (let's recall that the exact differential 
structure of such a moment-angle manifold is conjectured in \cite{F-W}).
\end{remark}

\begin{remark}
  We can notice that both polytopes on the first picture are flag polytopes. So 
this construction produces nonrigid flag polytopes, as announced.
\end{remark}

  Let's now prove the theorem:
  
\begin{demonstration}
  Let $P_1 $, $P_2 $ two polytopes joint by a puzzle-move.

  For $i=1,2$, let $\tilde{Z}_{{\Delta }_i }$ the preimage of ${\Delta }_i $ in 
$Z_{P_i }$ by the natural map, $Q_i = P_i \cap H_+ $, $R_i = P_i \cap H_- $.

  The two preimages $\tilde{Z}_{Q_1 }$ and $\tilde{Z}_{Q_2 }$ are 
diffeomorphic, as well as $\tilde{Z}_{R_1 }$ and $\tilde{Z}_{R_2 }$, 
$\tilde{Z}_{Q_i }$ and $\tilde{Z}_{R_i }$ having $\tilde{Z}_{{\Delta }_i }$ as 
common boundary.

  Now, we identify $Z_{Q_i }$ (resp. $Z_{R_i }$) to a manifold $\tilde{Z}_Q $ 
(resp to a manifold $\tilde{Z}_R $), yielding identifications $a_{i,+}$ (resp. 
$a_{i,-}$ of $\Delta _i $ with $\Delta $. With respect to the terminology 
thereabove, we have $i_+ = (a_{2,+})^{-1} \circ a_{1,+}$ whereas 
$i_- = (a_{2,-})^{-1} \circ a_{1,-}$.

  Both $Z_{P_1 }$ and $Z_{P_2 }$ are obtained by glueing $\tilde{Z}_Q $ and 
$\tilde{Z}_R $ along their common boundary $\tilde{Z}_{\Delta }$. Indeed, 
$\tilde{Z}_{\Delta }$ is given by the product of $Z_{\Delta }$ by a real torus 
$T^{\cal S }$, where $S$ corresponds to the set of facets of $Q$ and $R$ that 
are disjoint from $\Delta $.

  We can find Collar neigbourhoods 
$N_Q \simeq \tilde{Z}_{\Delta } \times T^{\cal S } \times [0;1[$ in 
$\tilde{Z}_Q $ and 
$N_R \simeq \tilde{Z}_{\Delta } \times T^{\cal S } \times [0;1[$ in 
$\tilde{Z}_R $ so that $Z_{P_i }$ is given by identifying (glueing) their 
interiors in the following way:
$$\begin{array}{ccc}
\tilde{Z}_{\Delta } \times T^{\cal S } \times [0;1[ & : &
\tilde{Z}_{\Delta } \times T^{\cal S } \times [0;1[ \\
(x , \gamma , t) & \sim & (\psi _i (x) , \gamma , 1-t)
\end{array}$$
where $\psi _i $ is the diffeomorphism of $Z_{\Delta }$ induced by 
the automorphism $a_i = a_{i,-} \circ (a_{i,+})^{-1}$ of $\Delta $ associated 
to $P_i $. As by hypothesis $(i_- )^{-1} \circ i_+ $ is harmless, the 
composition $a_1 \circ (a_2 )^{-1}$ is harmless too, so, by 
proposition~\ref{caract-harm}, $\psi_1 \circ (\psi_2 )^{-1}$ is isotopic to 
identity, in other words  $\psi_1 $ and $\psi_2 $ are isotopic. Hence the 
manifolds $Z_{P_1 }$ and $Z_{P_2 }$ are actually diffeomorphic.

  Moreover, this diffeomorphism is graded.

  Indeed $\tilde{Z}_{{\Delta }_i }$ cuts $Z_{P_1 }$ into two parts. The 
homology of $Z_{P_1 }$ is then, by Mayer-Vietoris theory, generated by two 
kinds of classes:
\begin{itemize}
\item Classes induced by one of the two parts. In this case, this class is 
induced by a class that has a representant in $\partial P_1 $ in one side of 
the hyperplane containing $\Delta $. This class is identified with a class of 
$P_2 $ on the same side (we consider the puzzle-move fixes this side and moves 
the other). Obviously, these two classes have equal degrees.

\item Classes (pure) yielding nontrivial classes of $\tilde{Z}_{{\Delta }_1 }$. 
Such a class and its correspondant by the diffeomorphism yield "the same" class 
of $\tilde{Z}_{\Delta }$. As the degree of a representant of this class of 
$\tilde{Z}_{\Delta }$ is one less than the degree of a representant of the 
class itself, the two classes have equal degrees.
\end{itemize}

  This proves that the diffeomorphism actually preserves bigraduation.
\end{demonstration}

  As a first application, we can again recover that two polytopes obtained by 
cutting a vertex to the same one are $Gr$-equivalent, as we can see on the 
following picture which is valid in any dimension~:

\scalebox{.7}{\includegraphics{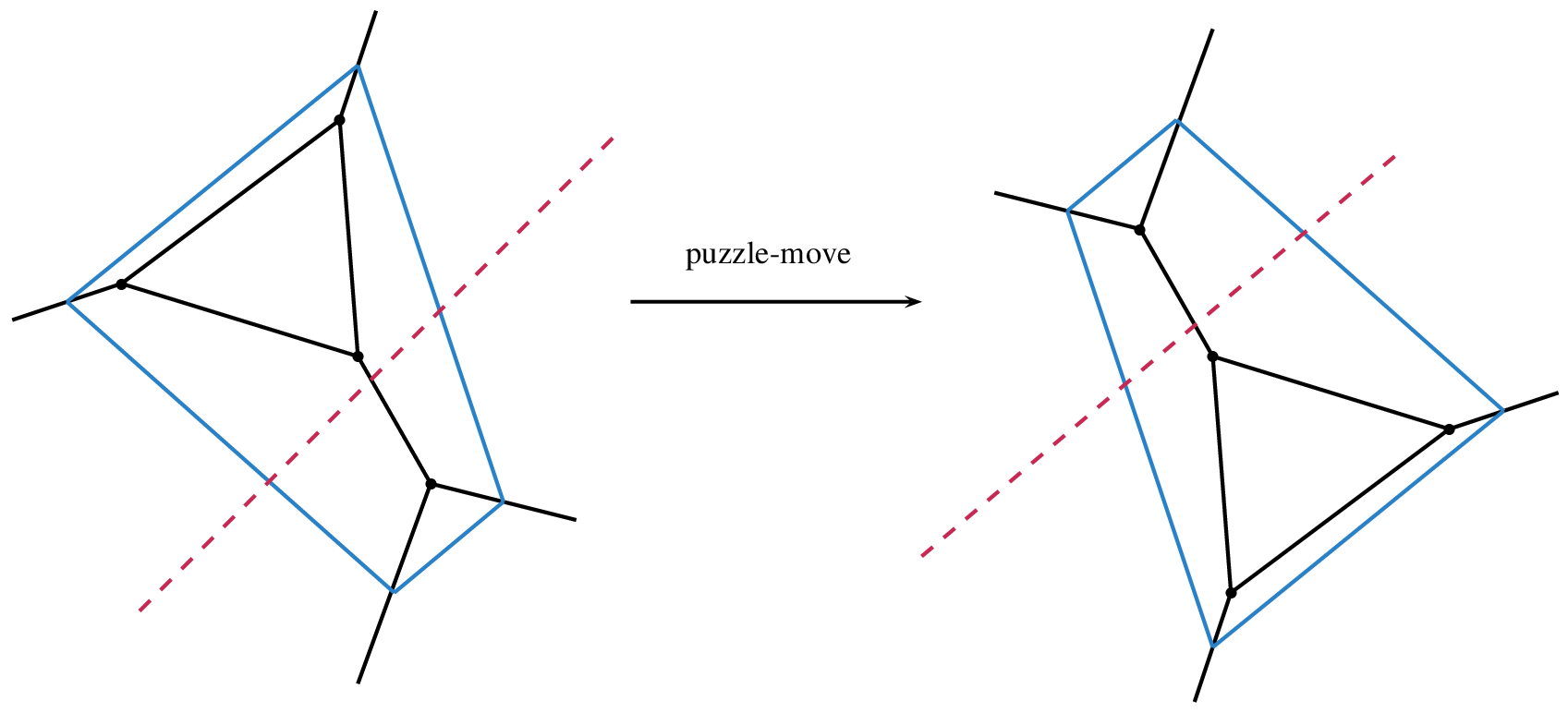}}

\paragraph{Dimension 3}

  We focus here on the case of puzzle-moves of $3$-dimensional polytopes.

  If $P$ is a $3$-dimensional polytope, then its intersection with a hyperplane 
is a polygon. The only nontrivial harmless automorphisms of polygons are 
obtained for a triangle (any permutation of sides is suitable) or for a 
quadrilateral (each side must then be sent either on itself or on its 
opposite), so a nontrivial puzzle-move from $P$ must benefit from such 
configurations.

  The case of a triangle implies that $P$ is decomposable (as a connected sum) 
and the permutation describes the difference between the glueings (at the same 
vertices).

  The case of a quadrilateral is perhaps less easily described, though quite 
clear on a picture. Notice it requires a $4$-belt on $P$, that separates it 
into two polytopes that are not too trivial (for example, a $4$-belt 
surrounding just a quadrilateral facet of $P$ produces a trivial puzzle-move).

  We know that polyhedra without $3$- nor $4$-belts (so-called Pogorelov 
polyhedra) are $Gr$-rigid \cite{F-M-W}, so $Gr$-equivalence between two 
polyhedra requires such kind of belts. Moreover, all known constructions of 
$Gr$-equivalent polytopes (connected sums, biflips, ...) are, in dimension $3$, 
obtained by puzzle-moves.

  Indeed, in dimension $3$, a biflip of edges can be obtained by two 
puzzle-moves like on the following picture:

\scalebox{.75}{\includegraphics{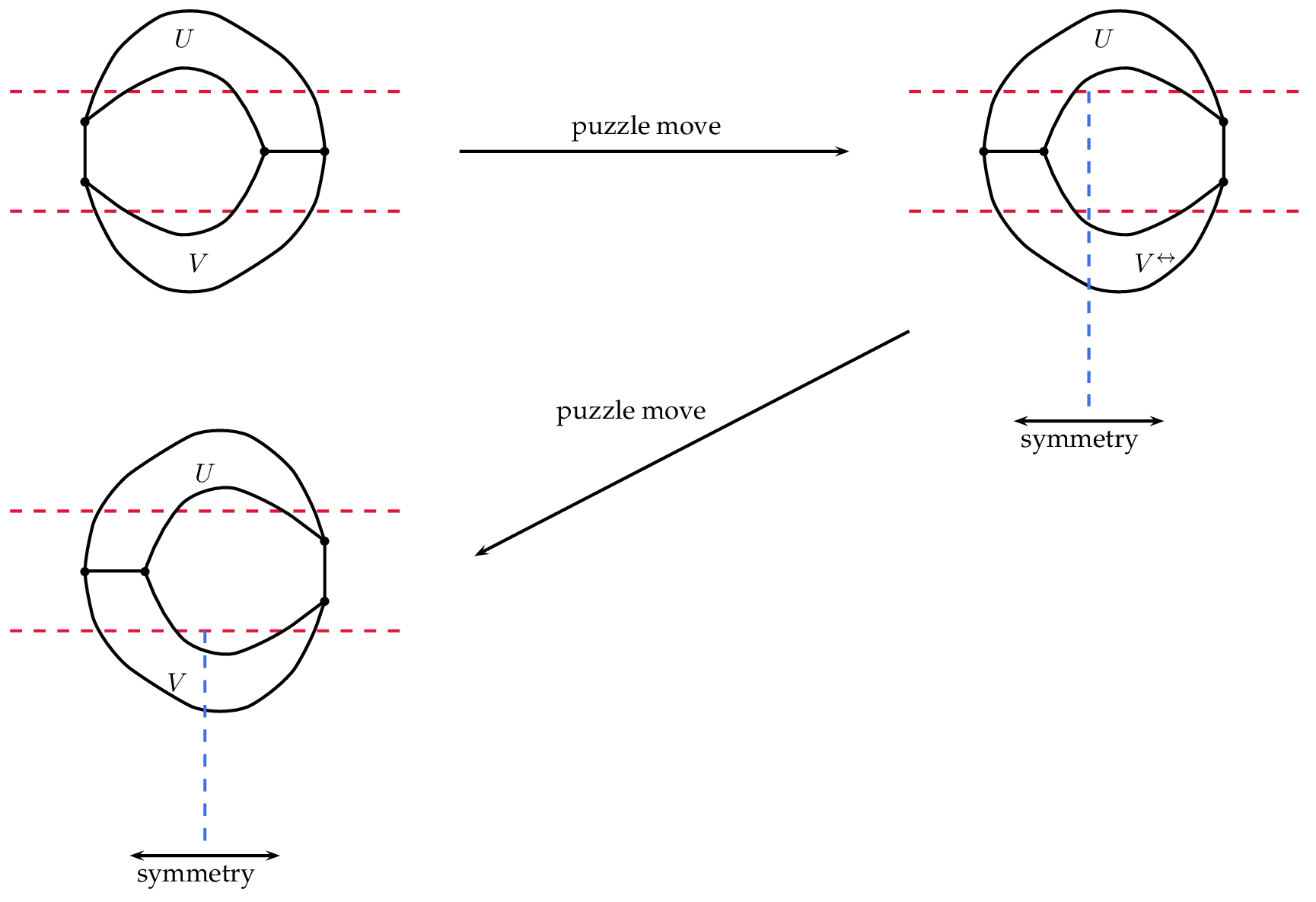}}

  This leads to conjecture:

\begin{conjecture}
  For $3$-dimensional polytopes, $Gr$-equivalence coincides with 
puzzle-equivalence.
\end{conjecture}

\subsection{Puzzle-rigidity}

\begin{definition}
  We say that a polytope is puzzle-rigid if it is puzzle-equivalent to no other 
polytope.
\end{definition}

  Naturally, puzzle-rigidity is weaker than $Gr$-rigidity.

  The following lemma may provide an easy way to gaurantee puzzle-rigidity:

\begin{lemma}
\label{facettes-pour-puzzle}
  Let $P$ a polytope, $F$ and $G$ two facets of $P$. Assume that the subgraph 
of the $1$-skeleton of $P$ induced by the vertices that are neither on $F$ nor 
on $G$ is connected.

  Then there is no nontrivial puzzle-move induced by the transposition $(F,G)$.

  If any pair of facets of $P$ satisfies the condition thereup, the polytope 
$P$ will be called {\em{$W$-puzzle-rigid}} and, in particular, it is 
puzzle-rigid. 
\end{lemma}

\begin{proof}
  Consider a puzzle-move induced by a transposition $(F,G)$. If $F$ and $G$ 
satisfy the condition of the lemma, then all the vertices outside $F \cup G$ 
lie on the same side of the separator as no vertex joining two such vertices 
can cross it. This side can be chosen as the right one. On the left side, every
vertex lies on $F$ or $G$ and if such a vertex $v$ is not on $G$, it is
adjacent to exactly one vertex $v'$ that is not on $F$. As the edge joining 
$v$ to $v'$ does not lie on $F$ nor $G$, it cannot cross $H$, so $v'$ must be 
on the left side and as it is not on $F$, it must lie on $G$.

  This proves that the transposition $(F,G)$ preserves the left vertices, so it 
produces a trivial puzzle-move.

  The second part of the lemma follows immediately.
\end{proof}

  We now use this lemma to highlight the difference between $Gr$-equivalence 
and puzzle-equivalence.
  
  In \cite{Bo2}, we classified $Gr$-rigidity among polytopes with three facets
more than their dimension. Recall that generically, such a polytope is 
nonrigid. Hence, the following theorem shows a neat difference between 
$Gr$-rigidity and puzzle-rigidity.

\begin{proposition}
  Any polytope with (at most) three facets more than its dimension is 
puzzle-rigid.
\end{proposition}

\begin{proof}
  We just have to verify $W$-puzzle-rigidity, i.e. the criterion of 
lemma~\ref{facettes-pour-puzzle} for such polytopes.

  If $n = d+1$, then any vertex lies on every facet but one, so in the union 
of any pair of facets.

  If $n = d+2$, then any vertex lies on every facet but two, so there is at 
most one vertex not lying on the union of a given pair of facets.

  If $n = d+3$, then any vertex of the complement of the union of two given 
facets lies on exacly all other facets but one, so any two such vertices are 
adjacent, i.e. the corresponding graph is complete, which proves the criterion.
\end{proof}

\paragraph{Puzzle-rigidity and wedges}

  We now deal with relations between puzzle-rigidity, or $W$-puzzle-rigidity, 
and the wedging operation.

  For the second one, the situation is particularly simple:

\begin{proposition}
  Let $P$ a polytope, $F$ a facet of $P$. Then $W_F P$ is $W$-puzzle-rigid if 
and only if $P$ is $W$-puzzle-rigid.
\end{proposition}

\begin{proof}
  Assume $P$ is $W$-puzzle-rigid and $F$ is a facet of $P$. Consider now two 
distinct facets $\tilde{G}$ and $\tilde{H}$ of $W_F P$ corresponding to facets 
$G$ and $H$ of $P$. If $G=H=F$, then $\tilde{G}$ and $\tilde{H}$ are the main 
facets, so their union contains all vertices of $W_F P$. If $G=F$, $H \neq F$, 
then the graph on the complement of $\tilde{G} \cup \tilde{H}$ is isomorphic to 
the graph on the complement of $F \cup H$, which is connected by hypothesis.

   If neither$G$ nor $H$ is $F$, then both intersections of the graph induced 
by the complement of $\tilde{G} \cup \tilde{H}$ with the main facets are 
connected and any vertex of any of these graphs is adjacent to one of the other 
graph (by swapping the main facets). Anyway, the desired graph is actually 
connected.

  Conversely, assume $P$ is not $W$-puzzle-rigid and let $G,H$ two facets of 
$P$ so that the graph induced by their complement is disconnected.

  Consider now facets $\tilde{G}$, $\tilde{H}$ over respectively $G$ and $H$. 
If $G$ is $F$, then the graph induced by the complement of 
$\tilde{G} \cup \tilde{H}$ is isomorphic to this disconnected graph. If neither
$G$ nor $H$ is $F$, then, as two adjacent vertices of wedge have their 
projections on the basic polytope equal or adjacent, points projecting on 
different components of the (graph induced by the) complement of $G \cup H$ 
cannot be adjacent, so have to be in different components of the (graph induced 
by the) complement of $\tilde{G} \cup \tilde{H}$. Hence, the desired graph is 
disconnected.
\end{proof}

  As we shall see, the situation is not that simple concerning puzzle-rigidity. 
For instance, an hexagon is puzzle-rigid, but an hexagonal book (a wedge over 
the hexagon) is not. We can ask:

\begin{question}
  Is $W$-puzzle-rigidity of a polytope equivalent to puzzle-rigidity of all 
multiwedges over it?
\end{question}

  In fact, we can see that a puzzle-move on a polytope induces naturally a 
combinatorial operation on any (multi)wedge over it.

\begin{definition}
\label{moveinduit}
  Let $(P , Delta ^+ , \phi )$ a puzzle-move, $F$ a facet of $P$ and $H$ the 
space of $\Delta $. We can consider $\bar{F} = F \cap H$. There is a hyperplane 
$H'$ of the space of $W_F P$ that separates vertices according to whether their 
projection on $P$ is in $H_+ $ or $H_- $, and let's make the natural choice for 
$H'_+ $ and $H'_- $. Combinatorially, its intersection with $W_F P$ is 
$W_{\bar{F}} (P \cap H)$. We now consider an automorphism $\phi '$ of 
$W_{\bar{F}} (P \cap H)$ that fixes one of the two main facets and acts like 
$\phi $ on the others (throw the natural identification).

  Then $\phi '$ is a harmless automorphism of $W_{\bar{F}} (P \cap H)$, and we 
can consider the simplicial complex on the facets of $W_F P$ whose maximal 
elements correspond to either the facets contaning a vertex in $H'_+ $ or those 
containing a vertex in $H'_- $ but where the ones meeting $H'$ are repaced by 
their image by $\phi '$ (we just mimic the combinatorial realisation of a 
puzzle-move).

  This complex will be noted $P_{H'^+ , \phi '}$.
\end{definition}

  In many cases, we can guarantee that the passage from $P$ to 
$P_{H'^+ , \phi '}$ is a nontrivial puzzle-move, so it prevents the rigidity of 
the wedge, and sometimes of any multiwedge over the base polytope.

  The following proposition often applies when dealing with (multi)wedges:

\begin{lemma}
\label{wedge-non-rig}
  Let $P$ a polytope. Assume that there is a puzzle-move that transforms $P$ 
into a polytope with more wedge-equivalence classes. Then the operation induced 
by proposition~\ref{moveinduit} on a (multi)wedge over $P$ cannot be a trivial 
puzzle-move.
\end{lemma}

\begin{proof}
  Indeed, the number of wedge-classes of $W_{\alpha } P$ is equal to the one of 
$P$ and, calling $P_{\Delta }$ the polytope obtained from $P$ by the required 
puzzle-move, the object obtained by proposition~\ref{moveinduit} on 
$W_{\alpha } P$, if a polytope, contains $P_{\Delta }$ as face, so has at least 
asmany wedge-classes as it, i.e. more than $W_{\alpha } P$. So we cannot get 
isomorphic polytopes.
\end{proof}

  Realisability of moves may be more delicate, but we can settle it in two 
relevant cases:

\begin{lemma}
  If $F$ is disjoint from $H$ or if $\bar{F}$ is fixed by $\phi $, then the 
operation is a puzzle-move.

  If the automorphism $\phi $ of $\Delta $ can be extended to an isometry of a 
neigbourhood of $\Delta $ in $P$, then the operation is a puzzle-move.
\end{lemma}

\begin{remark}
  If the original puzzle-move consists of changing the glueing of a 
connected sum of two polytopes (in other words if $\Delta $ is a simplex), 
then there is realisation of $P$ so that a neigbourhood of $\Delta $ is the 
product of a reguular simplex and an interval. In this case the proposition 
thereup applies so the operation of proposition~\ref{moveinduit} is actually a 
puzzle-move.
\end{remark}

\begin{proof}
  In the firts case, we easily notice that the produced object is simply 
$W_F Q$. So it is a polytope.

  In the second case, the isometry of the neigbourhood of $\Delta $ in $P$ 
naturally extends to an isometry of a neigbourhood of $H' \cap W_F P'$ in $P'$, 
and the corresponding glueing also produces a convex body, i.e. a polytope.
\end{proof}

  As an application, we have:

\begin{theorem}
  Let $n \geq 6$. Then no nontrivial multiwedge over the 
$n$-gon is puzzle-rigid.
\end{theorem}

\begin{demonstration}
  Indeed, if we consider a simple wedge over an hexagon, i.e. a hexagonal book, 
we can notice that a puzzle-move corresponding to a connected sum can turn a 
wedge over an hexagon into a polytope which is not a nontrivial wedge, so has 
one more wedge-class.

\scalebox{.75}{\includegraphics{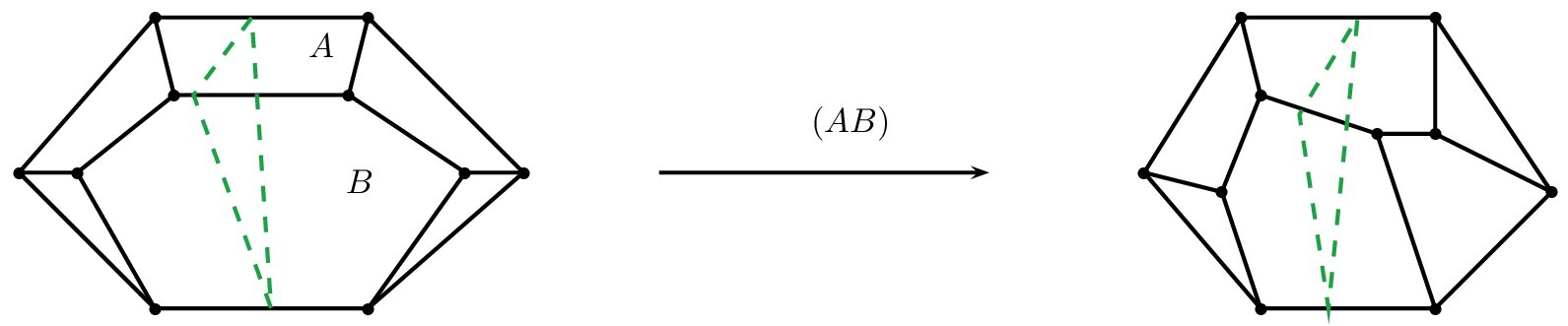}}

The lemmas thereup permit to conclude.
\end{demonstration}

Note: We can see that dual cyclic polytopes of dimension $d \geq 3$ and at 
least $d+4$ facets are not puzzle-rigid, hence not rigid, neither are any 
(multi)wedge over such a polytope. The same question is open by replacing 
"dual cyclic" by "dual neighbourly".

\paragraph{Puzzle-rigidity and products}

  The puzzle-equivalence between products of polytopes is the simplest we can 
expect. Precisely, we get:

\begin{proposition}
  Let $P$ and $Q$ two polytopes. Then any polytope puzzle-equivalent to 
$P \times Q$ is a product $P' \times Q'$ where $P'$ is puzzle-equivalent to 
$P$ and $Q'$ puzzle-equivalent to $Q$.
\end{proposition}

\begin{proof}
Let $P$ and $Q$ two polytopes, $F$ a facet of $P$ and $G$ a facet of $Q$. Then
$F$ and $G$ satisfy the condition of lemma~\ref{facettes-pour-puzzle}. Indeed, 
if we consider two vertices $(v_{P,1} , v_{Q_1 })$ and $(v_{P,2} , v_{Q,2})$
out of $F \cup G$, it means that neither $v_{P,i}$ is on $F$ and neither
$v_{Q,i}$ is on $G$. Then there is a path in $P$ from $v_{P,1}$ to $v_{P,2}$ 
avoiding $F$ and a path in $Q$ from $v_{Q,1}$ to $v_{Q,2}$ avoiding $G$. Using 
these paths, we get a path from $(v_{P,1} , v_{Q_1 })$ to $(v_{P,2} , v_{Q,2})$ 
avoiding $F \cup G$.

  So any elementary puzzle-move of $P \times Q$ permutes two facets of the same 
kind ($F \times Q$ or $P \times G$), we can assume to be 
$\tilde{F} = F \times Q$ and $\tilde{F'} = F' \times Q$. Let's then prove that
the polytope $R$ obtained after this puzzle-move has the form 
$\tilde{P} \times Q$.

  Consider a vertex $v$ of $P$. If it lies neither on $F$ nor on $F'$, then for 
any $v_Q $ on $Q$, the vertex $(v , v_Q )$ is preserved by the move, so is a 
vertex of $R$. Another thing to notice is this case is that all these vertices
lie on the same side of $H$, otherwise a vextex joining two of them should
cross $H$, which is incompatible with the existence of the puzzle-move.

  If $v \in F \cap F'$, then for any $v_Q $ on $Q$, the vertex $(v , v_Q )$ is 
preserved by the move, whichever side of $H$ $v$ lies, so is also a vertex of 
$R$.

  If $v$ is on $F$, not on $F'$, then there is exactly one vertex $\tilde{v}$ 
of $P$ that is adjacent to $v$ and not on $F$. We just distinguish two cases: 

  If $\tilde{v}$ is on $F'$, then for any vertex $v_Q $ of $Q$, the pair of 
vertices $\{ (v,v_Q ) ; (\tilde{v}, v_Q ) \} $ is necessarily preserved by the 
puzzle-move, so both, and in paticular $(v,v_Q )$ are vertices of $R$.

  If now $\tilde{v}$ is not on $F'$, then, as we have mentioned, all vertices 
$(\tilde{v}, v_Q )$ are on the same side of $H$. As no edge joining 
$(v , v_Q )$ to $(v' , v_Q )$ can cross $H$ as it neither lies on $F$ nor $F'$, 
all the vertices of the form $(v , v_Q )$ must be on the same side of $H$. So
applying the puzzle-move does not break the product-by-$Q$ structure.

  The same argument also being valid for vertices lying on $F' - F$, the 
polytope obtained after the puzzle-move is actually a product of some polytope 
$\tilde{P}$ by $Q$.

  The end of the proof is easy, we can fix any particular vertex $v_Q $ of $Q$ 
and we observe that when we consider the only face $P \times \{ v_Q \} $ of 
$P \times Q$, the global puzzle-move induces a puzzle-move from $P$ to 
$\tilde{P}$ induced by the transposition $(F,F')$.

  The proposition is then proved.
\end{proof}

  In particular, the product of two puzzle-rigid polytopes is also puzzle-rigid.

  Notice that replacing puzzle-rigid or puzzle-equivalent by rigid or 
$Gr$-equivalent in the previous proposition leads to open questions.

\paragraph{Generalisation?}

  Puzzle-equivalence is defined for embedded polytopes by introducing an 
hyperplane of the ambiant space. Naturally, we can suspect that there exists 
and would be valuable to produce a more general, purely combinatorial version 
of this operation that would work on more general combinatorial objects than 
polytopes, and induce $Gr$-equivalence on corresponding moment-angle manifolds.

{\footnotesize {Bosio Fr\'ed\'eric \\
Universit\'e de Poitiers \\
UFR Sciences SP2MI \\
D\'epartement de Math\'ematiques \\
UMR CNRS 6086 \\
Teleport 2 \\
Boulevard Marie et Pierre Curie \\
BP 30179 \\
86962 Futuroscope Chasseneuil CEDEX 

e-mail~: bosio@math.univ-poitiers.fr}}

\end{document}